\documentclass[11pt]{amsart}

\usepackage{amscd}
\usepackage{amsmath, amssymb}
\usepackage{amsfonts}

\newcommand{\de}{\partial} 
\newcommand{\db}{\overline{\partial}}

\newcommand{\Ric}{\mathrm{Ric}}

\newcommand{\ov}[1]{\overline{#1}}

\newcommand{\mn}{\sqrt{-1}}

\newcommand{\vp}{\varphi}
\newcommand{\vol}{\mathrm{Vol}}
\newcommand{\diam}{\mathrm{diam}}

\renewcommand{\leq}{\leqslant}
\renewcommand{\geq}{\geqslant}

\begin{document}
\newcounter{remark}
\newcounter{theor}
\setcounter{remark}{0}
\setcounter{theor}{1}
\newtheorem{claim}{Claim}
\newtheorem{theorem}{Theorem}[section]
\newtheorem{proposition}{Proposition}[section]
\newtheorem{question}{question}[section]
\newtheorem{lemma}{Lemma}[section]
\newtheorem{defn}{Definition}[theor]

\newtheorem{corollary}{Corollary}[section]

\newenvironment{example}[1][Example]{\addtocounter{remark}{1} \begin{trivlist}
\item[\hskip
\labelsep {\bfseries #1  \thesection.\theremark}]}{\end{trivlist}}

\title[Degenerations of Calabi-Yau metrics]{Degenerations of Calabi-Yau metrics}
\author{Valentino Tosatti}
 \address{Department of Mathematics \\ Columbia University \\ New York, NY 10027}

  \email{tosatti@math.columbia.edu}
\begin{abstract} 
We survey our recent work on degenerations of Ricci-flat K\"ahler metrics on compact Calabi-Yau manifolds with K\"ahler classes approaching the boundary of the K\"ahler cone. 
\end{abstract}
\thanks{I am grateful to Yuguang Zhang for some very useful comments. This work was partially supported by National Science Foundation grant DMS-1005457.}

\maketitle

\section{Introduction}
A compact K\"ahler manifold $X$ of complex dimension $n$ is called a {\em Calabi-Yau manifold} if its first Chern class
$c_1(X)$ vanishes in the cohomology group $H^2(X,\mathbb{R})$. This is equivalent to requiring that the canonical 
bundle $K_X$ be torsion, so that $K_X^{\otimes \ell}\cong \mathcal{O}_X$ for some integer $\ell\geq 1$. The following are some
simple examples of Calabi-Yau manifolds.\\

\begin{example} Let $X=\mathbb{C}^n/\Lambda$ be the quotient of Euclidean space $\mathbb{C}^n$ by a lattice
$\Lambda\cong\mathbb{Z}^{2n}$. Then $X$ is topologically just a torus $(S^1)^{2n}$ and it has trivial tangent bundle
and therefore also trivial canonical bundle. All Calabi-Yau manifolds of complex dimension $n=1$ are tori, and
are also called {\em elliptic curves}.
\end{example}

\begin{example} A Calabi-Yau manifold with complex dimension $n=2$ which is also simply connected is called a 
{\em K3 surface}. Every Calabi-Yau surface is known to be either a torus, a $K3$ surface, or a finite unramified
quotient of these. In general these quotients will have torsion but nontrivial canonical bundle, as is the case for
example for {\em Enriques surfaces} which are $\mathbb{Z}/2$ quotients of $K3$.
\end{example}

\begin{example} Let $X$ be a smooth complex hypersurface of degree $n+2$ inside complex projective space $\mathbb{CP}^{n+1}$.
Then by the adjuction formula the canonical bundle of $X$ is trivial, and so $X$ is a Calabi-Yau manifold. When $n=1$
we get an elliptic curve and when $n=2$ a $K3$ surface. More generally one can consider smooth
complete intersections in product of projective spaces, with suitable degrees, and get more examples of Calabi-Yau manifolds. \end{example}

\begin{example} Let $T=\mathbb{C}^2/\Lambda$ be a torus of complex dimension $2$ and consider the reflection through
the origin $i:\mathbb{C}^2\to \mathbb{C}^2$. This descends to an involution of $T$ with $16$ fixed points (the $2$-torsion points of $T$),
and we can take the quotient $Y=T/i$ which is an algebraic variety with $16$ singular rational double points (also known as 
orbifold points). We resolve these $16$ points by blowing them all up and we get a map $f:X\to Y$ where $X$ is a smooth $K3$ surface,
known as the {\em Kummer surface} of the torus $T$.
\end{example}

The fundamental result about Calabi-Yau manifolds, which is also the reason for their name, is the following

\begin{theorem}[Yau's solution of the Calabi Conjecture, 1976 \cite{yaupnas, yau1}]\label{yauthm}
On any compact Calabi-Yau manifold $X$ there exist K\"ahler metrics with Ricci curvature identically zero,
which naturally have restricted holonomy contained in $SU(n)$. Moreover, there is a unique such Ricci-flat metric in each
K\"ahler class of $X$.
\end{theorem}

To make this statement more precise, recall that a K\"ahler metric $g$ together with the complex structure $J$ of $X$ defines 
a real closed $2$-form $\omega$ (the K\"ahler form) by the formula $\omega(X,Y)=g(JX,Y)$. Conversely we can recover $g$ from $\omega$
by the formula $g(X,Y)=\omega(X,JY)$, and so we will often refer to $\omega$ simply as the K\"ahler metric.
The form $\omega$ is also $J$-invariant, in the 
sense that $\omega(JX,JY)=\omega(X,Y)$, which also means that it is of complex type $(1,1)$. Therefore it defines a cohomology
class
$$[\omega]\in H^2(X,\mathbb{R})\cap H^{1,1}_{\ov{\de}}(X)=:H^{1,1}(X,\mathbb{R}).$$
A cohomology class $\alpha$ in $H^{1,1}(X,\mathbb{R})$ which can be written as $\alpha=[\omega]$ for some K\"ahler metric $\omega$ is
called a {\em K\"ahler class}. The set of all K\"ahler classes is called the {\em K\"ahler cone} of $X$ and is
an open convex cone $$\mathcal{K}_X\subset H^{1,1}(X,\mathbb{R}),$$
which has the origin as its vertex. 

With these notions in place, we can restate Theorem \ref{yauthm} by saying that on a compact Calabi-Yau there is 
a unique Ricci-flat K\"ahler metric in each K\"ahler class $\alpha\in\mathcal{K}_X$. These metrics are almost never explicit, and Yau constructed them by solving a fully nonlinear complex Monge-Amp\`ere PDE.
Since the K\"ahler cone $\mathcal{K}_X$
is open, the following question is very natural\\

\noindent {\bf Question 1.1 }What is the behaviour of these Ricci-flat K\"ahler metrics when the class $\alpha$
degenerates to the boundary of the K\"ahler cone?\\

This question was posed by many people, including Yau \cite{yau2, yau3}, Wilson \cite{wilson} and McMullen \cite{ctm}. To get a feeling for what the K\"ahler cone and its boundary represent geometrically, we start with the following observation.
If $V\subset X$ is a complex subvariety of complex dimension $k>0$, then 
it is well known (from the work of Lelong) that $V$ defines a homology class 
$[V]$ in $H_{2k}(M,\mathbb{R})$. Moreover if $[\omega]$ is a K\"ahler class, the pairing $\langle [V],[\omega]^{\smile k}\rangle$ equals
$$\int_V \omega^k=\vol(V,\omega)>0,$$
the volume of $V$ with respect to the K\"ahler metric $\omega$ (Wirtinger's Theorem). It follows that if a class $\alpha$ is on the boundary of $\mathcal{K}_X$ and if $V$ is any complex subvariety then the pairing $\langle [V],\alpha^{\smile k}\rangle$ is nonnegative, and moreover a theorem of Demailly-P\u{a}un \cite{dp} shows
that there must be subvarities $V$ with pairing zero. Therefore as we approach the class $\alpha$ from inside $\mathcal{K}_X$, these subvarieties have volume that
goes to zero, and the Ricci-flat metrics must degenerate (in some way) along these subvarieties.\\

We now make Question 1.1 more precise. On a compact Calabi-Yau manifold $X$
fix a nonzero class $\alpha_0$ on the boundary of $\mathcal{K}_X$ and let $\{\alpha_t\}_{0\leq t\leq 1}$ be a smooth path of classes in $H^{1,1}(X,\mathbb{R})$ originating at $\alpha_0$ and with $\alpha_t\in \mathcal{K}_X$ for $t>0$. Call $\omega_t$ the unique Ricci-flat K\"ahler metric on $X$ cohomologous to $\alpha_t$ for $t>0$, which is produced by Theorem \ref{yauthm}.\\

\noindent {\bf Question 1.2 }What is the behaviour of the Ricci-flat metrics $\omega_t$ when $t$ goes to zero?\\

Of course we could also consider sequences of classes instead of a path, and all we are going to say in this paper works equally well in that case.
Notice that we are not allowing the class $\alpha_t$ to go to infinity in $H^{1,1}(X,\mathbb{R})$ as it approaches $\de\mathcal{K}_X$. Because of this, we can prove the following basic fact, independently discovered by Zhang \cite{zh}:

\begin{theorem}[T. \cite{deg}, Zhang \cite{zh}] The diameter of the metrics $\omega_t$ has a uniform upper bound as $t$ approaches zero,
\begin{equation}\label{diameter}
\diam(X,\omega_t)\leq C.
\end{equation}
\end{theorem}

On the other hand it is easy to construct examples of Ricci-flat K\"ahler metrics with unbounded cohomology class that violate \eqref{diameter}, by just rescaling a fixed metric by a large number. There are also examples where the class approaches the boundary of the K\"ahler cone: on the torus $\mathbb{C}^2/\Lambda=T^2\times T^2$ take the flat metric that gives one $T^2$ factor area $t$ and the other $T^2$ factor area $t^{-1}$.

Going back to Question 1.2, the problem splits naturally into two cases which exhibit a rather different behaviour, according to whether the total integral $\int_X \alpha_0^n$ is strictly positive or zero. If $\int_X \alpha_0^n$ is positive this means that the volume
$$\vol(X,\omega_t)=\int_X\omega_t^n=\int_X\alpha_t^n$$
remains bounded away from zero as $t\to 0$, and this is called the {\em non-collapsing} case. If $\int_X \alpha_0^n=0$ then the volume $\vol(X,\omega_t)$
converges to zero, and this is called the {\em collapsing} case.

The main question 1.2 falls into the general problem of understanding limits of sequences of Einstein manifolds with an upper bound for the diameter (but no bound for the sectional curvature in general), a topic that has been extensively studied (see e.g. \cite{gr, and, bkn, tian, CC}). Our results are of a quite different nature from these works, because the convergence that we get is in a stronger sense, we have uniqueness of the limit, and we do not need to modify the metrics by diffeomorphisms. On the other hand these other works apply in much more general setups, and are especially effective in complex dimension $n=2$.

One final comment: in our work we always fix the complex structure and vary the K\"ahler class. If instead one varies the complex structure as well the behaviour is expected to be much more complicated, except in complex dimension $n=2$ where changing the complex or K\"ahler structure are comparable operations because of an underlying hyperk\"ahler structure. In certain higher-dimensional cases, some convergence results have recently been obtained by Ruan-Zhang \cite{rz}.

\section{examples}\label{ex}
\setcounter{remark}{0}
First of all notice that Question 1.2 is only interesting if 
$$\dim H^{1,1}(X,\mathbb{R})>1,$$ because otherwise $\mathcal{K}_X$ reduces to
an open half-line and there is only one Ricci-flat K\"ahler metric on $X$ up
to global scaling by a constant, so the only possible degenerations are given by scaling this metric to zero or infinity. For this reason, the Question 1.2 is essentially void on Calabi-Yau manifolds of dimension $n=1$ (i.e. elliptic curves).\\

\begin{example} Let $X=\mathbb{C}^n/\Lambda$ be a complex torus. A Ricci-flat K\"ahler metric on $X$ is the same as a flat K\"ahler metric, and each flat metric can be identified simply with a positive definite Hermitian $n\times n$ matrix. The boundary of the K\"ahler cone is then represented by non-negative definite Hermitian matrices $H$ with nontrivial kernel $\Sigma\subset\mathbb{C}^n$ (notice that in this case every class on $\de \mathcal{K}_X$ has zero integral, so we are always in the collapsing case).

If the class $\alpha_0$ corresponds to such a matrix $H$ with the kernel $\Sigma$ which is $\mathbb{Q}$-defined modulo $\Lambda$, then we can quotient $\Sigma$ out and get a map $f:X\to Y=\mathbb{C}^m/\Lambda'$ to a lower-dimensional torus ($m<n$) such that $H=f^*H'$ with $H'$ a positive definite $m\times m$ Hermitian matrix. It follows that when $t$ approaches zero, the (Ricci-)flat metrics $\omega_t$ collapse to the flat metric on $Y$ that corresponds to $H'$.
Here collapsing has the precise meaning that the geometric limit (i.e. Gromov-Hausdorff limit) of $(X,\omega_t)$ has dimension strictly less than $n$.

If on the other hand the kernel $\Sigma$ is not $\mathbb{Q}$-defined, then $\Sigma$ defines a foliation on $X$ (which is not a fibration anymore) and the limit $H$ of the (Ricci-)flat metrics is a smooth nonnegative form which is {\em transversal} to the foliation (that means, positive in the complementary directions).\end{example}

\begin{example} Let $f:X\to Y$ be the Kummer $K3$ surface of a torus $T$,
where $Y=T/i$ is the singular quotient of $T$ and $f$ is the blowup map. Take $\alpha_0$ to be the pullback of an ample divisor on $Y$, and note that
$\int_X\alpha_0^2>0$. If we call $E$ the union of the $16$ exceptional divisors of $f$, that is the union of the $16$ spheres $S^2$ which are the preimages of the singular points of $Y$, then $E$ is a complex submanifold of $X$. Then Kobayashi-Todorov \cite{kt} proved that for any path $\alpha_t$ of K\"ahler classes that approach $\alpha_0$, the Ricci-flat metrics $\omega_t$ converge smoothly away from $E$ to the pullback of 
the unique flat orbifold metric on $Y$ cohomologous to the ample divisor we chose. Here an orbifold flat metric on $Y$ simply means a flat metric on $T$ which is invariant under $i$. Note that since the limit $Y$ has the same dimension as $X$, the Ricci-flat metrics are non-collapsing. This convergence result is proved using classical results on the moduli space of $K3$ surfaces, such as the Torelli theorem.
\end{example}

\begin{example} Let $X$ be a $K3$ surface which admits an elliptic fibration $f:X\to\mathbb{CP}^1=Y$. This means that $f$ is a surjective holomorphic map with all the fibers smooth elliptic curves except a finite number of fibers which are singular elliptic curves. Again we take $\alpha_0$ to be the pullback of an ample divisor on $Y$ and note that
$\int_X\alpha_0^2=0$. We also fix $[\omega]$ a K\"ahler class on $X$ and consider only paths which are straight lines of the form
$$\alpha_t=\alpha_0+t[\omega],$$
with $0\leq t\leq 1$. Again we call $\omega_t$ the unique Ricci-flat K\"ahler metric in the class $\alpha_t$ for $t>0$, and we call $E$ the union of all the singular fibers of $f$. Then Gross-Wilson \cite{gw} have shown that when $t$ goes to zero the metrics $\omega_t$ converge smoothly away from $E$ to the pullback $f^*\eta$, where $\eta$ is a K\"ahler metric on $Y=\mathbb{CP}^1$ minus the finitely many points $f(E)$ with singular preimage. Moreover they show that away from $E$ as $t\to 0$ we have
$$\omega_t\sim f^*\eta + t\omega_{SF}+o(t),$$
where $\omega_{SF}$ is a {\em semi-flat} form, that is a $(1,1)$-form that restricts
to a flat metric on each smooth torus fiber. More recently Song-Tian \cite{st} have noticed that the metric $\eta$ on $\mathbb{CP}^1\backslash f(E)$ satisfies 
$$\Ric(\eta)=\omega_{WP},$$
where $\omega_{WP}$ is the pullback of the Weil-Petersson metric from the moduli space of elliptic curves via the map that to a point in $\mathbb{CP}^1\backslash f(E)$ associates the elliptic curve which lies above that point. The $(1,1)$-form $\omega_{WP}$ is smooth away from $f(E)$ and is non-negative definite. If the fibers are all isomorphic elliptic curves then $\omega_{WP}$ vanishes identically; in this case $X$ cannot be $K3$ but instead it is the torus of Example \ref{ex}.1, and 
$\eta$ is the (Ricci-)flat metric $H'$ there. In general $\omega_{WP}$ measures the variation of the complex structure of the fibers of $f$.
\end{example}

\begin{example} McMullen \cite{ctm} has constructed a nonalgebraic $K3$ surface $X$, an automorphism $F:X\to X$ with infinite order, a nonempty open set $U\subset X$ and a real number $\lambda>1$ with the following property. If we fix any Ricci-flat K\"ahler metric $\omega$ on $X$, and we consider the Ricci-flat metrics
$$\omega_n=\lambda^{-n} (F^n)^*\omega,$$
then the cohomology classes $[\omega_n]$ converge to a nonzero limit class 
$\alpha_0\in\de\mathcal{K}_X$ with $\int_X\alpha_0^2=0$, and the metrics $\omega_n$ converge smoothly
to zero on $U$. The set $U$ is a {\em Siegel disk} for the automorphism $F$, which means that $U$ is $F$-invariant and it is biholomorphic to a disk where $F$ is conjugate to an irrational rotation. The number $\lambda$ is the largest eigenvalue for the action of $F^*$ on $H^{1,1}(X,\mathbb{R})$, and the class $\alpha_0$ is an eigenvector of $F^*$ with eigenvalue $\lambda$.
\end{example}

\section{Main Theorems}
Let $X$ be a compact Calabi-Yau manifold and $\alpha_0\in\de\mathcal{K}_X$ with $\int_X\alpha_0^n>0$. Let $E$ be the union of all complex subvarieties where $\alpha_0$ integrates to zero ($E$ itself is a complex subvariety). 

\begin{theorem}[T. \cite{deg}]\label{thm1}
In this situation there exists a smooth Ricci-flat K\"ahler metric $\omega_0$ on $X\backslash E$ such that for any path $\alpha_t$ as before, the Ricci-flat metrics $\omega_t$ with $t\to 0$ converge to $\omega_0$ smoothly on $X\backslash E$. Moreover if $\alpha_0\in H^2(X,\mathbb{Q})$ then there exist a birational map $f:X\to Y$ with $Y$ a singular Calabi-Yau variety such that $\omega_0=f^*\omega$ and $\omega$ is a singular Ricci-flat K\"ahler metric on $Y$.
\end{theorem}

A singular Calabi-Yau variety can be defined in algebraic geometry as a normal variety $Y$ with at worst canonical singularities such that some multiple of the canonical divisor $K_Y$ is Cartier and it is trivial. A singular Ricci-flat K\"ahler metric on such a space can be defined as a weak solution of the complex Monge-Amp\`ere equation, and its existence was proved by Eyssidieux-Guedj-Zeriahi \cite{egz1}.
Strictly speaking Theorem \ref{thm1} is only stated in \cite{deg} for projective Calabi-Yau manifolds, but it is possible to extend the arguments there to the more general K\"ahler case as stated here by using the recent work of Boucksom-Eyssidieux-Guedj-Zeriahi \cite{begz}. \\

This gives a possible answer to Question 1.2 in the noncollapsing case when 
$\int_X\alpha_0^n>0$. We now consider the collapsing case when $\int_X\alpha_0^n=0$. One major source of examples of such cohomology classes $\alpha_0$ is whenever we have a holomorphic fibration $f:X\to Y$ where $Y$ is a variety with lower dimension $m<n$, and we take $\alpha_0$ to be the pullback of an ample divisor on $Y$. Examples \ref{ex}.1 and \ref{ex}.3 above fall exactly in this category. A standard conjecture in algebraic geometry, the log abundance conjecture, would imply that whenever $\alpha_0\in \de\mathcal{K}_X$ with $\int_X \alpha_0^n=0$ satisfies $\alpha_0\in H^2(X,\mathbb{Q})$ then there is a fibration $f:X\to Y$ so that $\alpha_0$ is the pullback of an ample divisor on $Y$. So conjecturally this picture is the general picture for rational classes (compare also Example \ref{ex}.1, where rational classes give a fibration and irrational classes a foliation).

In this case we can always find a proper complex subvariety $E\subset X$ such that
$f:X\backslash E\to Y\backslash f(E)$ is a smooth submersion. The subvariety $E$ is given by the union of all singular fibers together
with all the fibers with dimensions strictly larger than $n-m$. This implies that for any $y\in Y\backslash f(E)$ the fiber $X_y=f^{-1}(y)$ is a smooth $(n-m)$-dimensional compact Calabi-Yau manifold. If we fix a K\"ahler metric $\omega$ on $X$ and use $\omega|_{X_y}$ as polarization, we get a map from $Y\backslash f(E)$ to the moduli space of polarized Calabi-Yau $(n-m)$-folds, analogously to Example \ref{ex}.3 above.

\begin{theorem}[T. \cite{deg2}]\label{thm2}
In this situation take $\alpha_t=\alpha_0+t[\omega]$. Then there exists a smooth K\"ahler metric $\eta$ on $Y\backslash f(E)$ such that the Ricci-flat metrics $\omega_t$ with $t\to 0$ converge to $f^*\eta$ on $X\backslash E$ in the $C^{1,\beta}$ topology of K\"ahler potentials (for any $0<\beta<1$). Moreover for any $y\in Y\backslash f(E)$ the metrics $\omega_t|_{X_y}$ converge to zero in the $C^1$ topology of metrics. The metric $\eta$ satisfies
$$\Ric(\eta)=\omega_{WP},$$
where $\omega_{WP}$ is the pullback of the Weil-Petersson metric from the moduli space of polarized Calabi-Yau $(n-m)$-folds.
\end{theorem}

The Weil-Petersson metric has the same properties that we discussed in Example \ref{ex}.3. One can give a more explicit formula for $\omega_{WP}$ as follows. Since the smooth fibers $X_y$ are Calabi-Yaus, there is an integer $\ell$ so that $K_{X_y}^{\otimes\ell}\cong \mathcal{O}_{X_y}$ for all $y\in Y\backslash f(E)$.
Thus we can find a holomorphic family (parametrized by $y\in Y\backslash f(E)$) of never-vanishing holomorphic sections $\Omega_y$ of $K_{X_y}^{\otimes\ell}$, and we get a volume form $(\Omega_y\wedge\ov{\Omega_y})^{1/\ell}$ on $X_y$. Then on $Y\backslash f(E)$ we have
$$\omega_{WP}=-\mn\de\db\log\int_{X_y}(\Omega_y\wedge\ov{\Omega_y})^{1/\ell}.$$
We note here that in the proof of Theorem \ref{thm2} the assumption that $\alpha_t=\alpha_0+t[\omega]$, which essentially means that $\alpha_t$ does not approach $\de\mathcal{K}_X$ tangentially, is used crucially in deriving the estimates which are described below.

\section{Discussion}
Theorems \ref{thm1} and \ref{thm2} are proved by showing suitable {\em a priori} estimates for a degenerating family of complex Monge-Amp\`ere equations. To be more precise, Yau's Theorem \ref{yauthm} is proved by solving the complex Monge-Amp\`ere equation of the form
$$(\omega+\mn\de\db\vp)^n=\Omega,$$
where $\omega$ is a fixed K\"ahler metric, $\Omega$ is a certain fixed smooth volume form and $\vp$ is a K\"ahler potential that we have to solve for, so the metric $\omega+\mn\de\db\vp$ is Ricci-flat. In this case Yau proved $C^3$ {\em a priori} estimates for the potential $\vp$, or equivalently $C^1$ estimates for the metric $\omega+\mn\de\db\vp$, and then deduced higher order estimates by a standard bootstrapping argument.

In the setting of Theorems \ref{thm1} and \ref{thm2}, we need to solve an equation
of the form
\begin{equation}\label{ma}
(f^*\omega_Y+t\omega_X+\mn\de\db\vp_t)^n=c_t\Omega,
\end{equation}
where $\omega_X, \omega_Y$ are fixed K\"ahler metrics on $X$ and $Y$ respectively,
$\Omega$ is a fixed smooth volume form on $X$ and $c_t$ is a constant that
is bounded away from zero in Theorem \ref{thm1} and is comparable to $t^{n-m}$ in Theorem \ref{thm2}.
Equation \eqref{ma} is equivalent to the fact that the K\"ahler metric $\omega_t=f^*\omega_Y+t\omega_X+\mn\de\db\vp_t$ is Ricci-flat. From now on we will focus on Theorem \ref{thm2}, where the analysis is much more complicated.
The equations \eqref{ma} are complex Monge-Amp\`ere equations that degenerate when $t$ approaches zero, in two different ways: first, the reference metrics
$f^*\omega_Y+t\omega_X$ degenerate, and second the right hand side approaches zero. In \cite{deg2} we first show that the K\"ahler potentials have uniformly bounded Laplacian on every compact set of $X\backslash E$ (making crucial use of estimates of Ko\l odziej \cite{kol} and extensions of these by Demailly-Pali \cite{DP} and Eyssidieux-Guedj-Zeriahi \cite{egz2}). Next, to prove collapsing, we show that 
the eigenvalues of the Ricci-flat metrics $\omega_t$
in the $n-m$ fiber directions are all of the order of $t$, so that the fibers are shrunk in the limit.  The remaining $m$ eigenvalues are of the order of $1$, so
the overall determinant of $\omega_t$ is exactly of
the order of $t^{n-m}$, as required by \eqref{ma}. More precisely, we prove that given any compact set $K\subset X\backslash E$ there is a constant $C_K$ so that on $K$ we have
$$C_K^{-1}(f^*\omega_Y+t\omega_X)\leq \omega_t\leq C_K (f^*\omega_Y+t\omega_X).$$
Moreover we show that when restricted to a fiber $X_y$ in $K$, the first derivatives of $\omega_t$ go to zero.
Once all the necessary {\em a priori} estimates are established, we can then pass to the limit weakly in \eqref{ma} and get exactly the equation for a metric on $Y\backslash f(E)$ with Ricci curvature equal to the Weil-Petersson metric. More details can be found in \cite{tesi, deg, deg2}.\\

Le us now mention a few open problems related to the above results, which seem very interesting and perhaps not too far from accessible.\\

\noindent{\bf Question 4.1 } It seems highly likely that if we consider the rescaled Ricci-flat metrics (the so-called {\em adiabatic limit}) $$\frac{\omega_t}{t}\bigg|_{X_y},$$
then these should converge to the unique Ricci-flat metric on $X_y$ in the cohomology class $\omega|_{X_y}$. If we denote by $\omega_{SF}$ the semi-flat form, which is a $(1,1)$-form on $X\backslash E$ that restricts
to the Ricci-flat metric on $X_y$ cohomologous to $\omega|_{X_y}$, then this would imply that as $t\to 0$
$$\omega_t\sim f^*\eta+t\omega_{SF}+o(t),$$
exactly as in Example \ref{ex}.3. Indeed if one could improve the convergence result in Theorem \ref{thm2}, say to convergence in the $C^{2}$ topology of K\"ahler potentials, then this would follow easily by taking the limit in the corresponding Monge-Amp\`ere equations. Unfortunately the convergence proved in Theorem \ref{thm2} does not seem to be strong enough to conclude this. The follwing is a stronger conjecture that as we said would directly imply Question 4.1:\\

\noindent{\bf Question 4.2 }In the setting of Theorem \ref{thm2} prove that the convergence of $\omega_t$ to $f^*\eta$ is actually in the smooth topology away from $E$. For this, it is enough to prove uniform $C^k$ estimates for $\omega_t$ on each compact set of $X\backslash E$, independent of $t>0$. In the proof of Theorem \ref{thm2} we have showed that we have uniform $C^0$ estimates, and $C^1$ in the fibers directions.\\

The natural remaining question is what happens to the Ricci-flat metrics $\omega_t$ when $\alpha_0$ is an irrational class with $\int_X\alpha_0^n=0$, so that
there is no fibration structure. We conjecture the following:\\

\noindent{\bf Question 4.3 }In this situation there is a proper complex subvariety $E\subset X$ and a smooth nonnegative
$(1,1)$-form $\omega_0$ on $X\backslash E$, which satisfies $\omega_0^n=0$, so that the Ricci-flat metrics $\omega_t$ converge smoothly away from $E$ to $\omega_0$.\\

In this case taking the kernel of $\omega_0$ we would get a foliation on $X\backslash E$ with leaves holomorphic subvarieties. In general this foliation will not be a holomorphic foliation, which means that the leaves will not vary holomorphically. In particular the dimension of the leaves
will not be constant, not even on a Zariski open set of $X$. One can see this in McMullen's example \ref{ex}.4, where (assuming that the metrics $\omega_n$ converge to a nonnegative form $\omega_0$) the foliation has $0$-dimensional leaves on the open set $U$, but it has nonzero dimensional leaves somewhere else, since $\alpha_0\neq 0$. Under the assumption that the sectional curvature of $\omega_t$ remains uniformly bounded, Ruan \cite{ruan} has shown that Question 4.2 is correct, and that moreover the foliation defined by $\omega_0$ is holomorphic. Therefore, McMullen's example shows that Ruan's result does not hold if the curvature is unbounded (the curvature of $\omega_n$ is of the order of $\lambda^n$).\\

One last problem that seems very interesting is whether the convergence in Theorems \ref{thm1} and \ref{thm2} holds in the Gromov-Hausdorff sense.
More precisely, in Theorem \ref{thm1} consider the metric space completion $(Z,d)$ of $(X\backslash E, \omega_0)$, while in Theorem \ref{thm2}
call $(Z,d)$ the metric space completion of $(Y\backslash f(E),\eta)$.\\

\noindent{\bf Question 4.4 }In the setting of either Theorem \ref{thm1} or \ref{thm2}, do the Ricci-flat manifolds
$(X,\omega_t)$ converge to $(Z,d)$ in the Gromov-Hausdorff sense? Moreover,
is $Z$ homeomorphic to $Y,$ the algebro-geometric limit?\\

The Gromov-Hausdorff converge is proved for $K3$ surfaces in the noncollapsing case in \cite{deg}, and there are further results in the noncollapsing case by Ruan-Zhang \cite{rz}. In the case of collapsing $K3$ surfaces, this was proved by Gross-Wilson \cite{gw}.

\end{document}